\documentclass[a4paper,leqno,11pt]{article}
\usepackage{amssymb}
\usepackage{amsmath}
\usepackage[german]{babel}
\begin{document}
\def\refname{ }
\begin{center}
\begin{Large}
\textbf{Algebraic independence of reciprocal sums}
\vskip0.1cm
\textbf{of certain Fibonacci-type numbers}
\end{Large}
\vskip0.7cm
\begin{large}
\large{Peter Bundschuh and Keijo V\"a\"an\"anen} 
\end{large}
\end{center}
\vskip0.7cm
\noindent
\textit{Abstract}: Let $(F_{n})_{n\geq 0}$ and $(L_{n})_{n\geq 0}$ denote the Fibonacci and Lucas 
sequences. The present work studies algebraic independence of the numbers 
$\sum_{h}b_h/F_{kr^h+\ell}, \sum_{h}c_h/L_{kr^h+\ell} \ (\ell\in\mathbb{Z},\ r\in\mathbb{N}\setminus\!\{1\})$, 
where $k\in\mathbb{N}$ is fixed, $(b_h)_{h\geq0},(c_h)_{h\geq0}$ are non-zero periodic sequences 
of algebraic numbers and the sum is taken over all integers $h \geq 0$ satisfying $kr^h+\ell>0$. 
Also results for more general sequences are obtained. The main tool is Mahler's method 
reducing the investigation of the algebraic independence of numbers (over $\mathbb{Q}$) to 
the one of functions (over the rational function field) if these satisfy certain types of functional 
equations. \\



\section{Introduction and main results}

The Fibonacci sequence $(F_{n})_{n\ge0}$ and the Lucas sequence $(L_{n})_{n\ge0}$ satisfy 
the recurrence
\[
y_{n+2}=y_{n}+y_{n+1}\ (n\geq 0)
\]
with initial conditions $y_{0}=0, y_{1}=1$ and $y_{0}=2, y_{1}=1$, respectively. In the present paper,
we are interested in the algebraic independence of the numbers
\begin{equation}\label{a1}
q_{\mu,r}:=\sum_{h=0}^\infty a_{h}\Big(\frac{1+\sqrt{5}}{2}\Big)^{-\mu kr^{h}},\ F_{\ell,r}:=
{\sum_{h=0}^\infty}\,\strut'\frac{b_{h}}{F_{kr^{h}+\ell}},\ L_{\ell,r}:={\sum_{h=0}^\infty}\,\strut'\frac{c_{h}}
{L_{kr^{h}+\ell}}
\end{equation}
\[
(\ell\in\mathbb{Z},\ r\in\mathbb{N}\setminus\!\{1\},\ \mu=1,\ldots,r-1) ,
\]
where $k\in\mathbb{N}$ is fixed, $(a_{h})_{h\ge0}, (b_{h})_{h\ge0}$ and $(c_{h})_{h\ge0}$ are 
non-zero periodic sequences of algebraic numbers, and $\sum_{h\ge0}^{\,'}$ is taken over those 
$h\in\mathbb{N}_{0}$ with non-zero denominator and $kr^{h}+\ell\ge0$. To present our first result,
we introduce a set of positive integers
\begin{equation}\label{a2}
\mathcal{D}=\{d\in\mathbb{N}\!:\ d\ne a^n \ (a,n\in\mathbb{N},\ n\ge2)\}.
\end{equation}
Then each integer $r\ge2$ has a unique representation in the form $r=d^{j}$ with $d\in\mathcal{D}, 
j\in\mathbb{N}$, and we define $d(r):=d$.\\

\textbf{Theorem 1}. \textit{The numbers
\begin{equation}\label{a3}
q_{\mu,r},\ F_{\ell,r},\ L_{\ell,r}\quad (\ell\in\mathbb{Z},\ r\in\mathbb{N}\setminus\!\{1\},\ \mu=1,\ldots,d(r)-1)
\end{equation}
are algebraically independent, unless $(b_{h})$ is constant. In this special case, $F_{0,2}$ is 
algebraic, but all other numbers in (\ref{a3}) are algebraically independent}.\\

This result gives a partial generalization of many earlier results considering separately 
the numbers $F_{\ell,r}$ or $L_{\ell,r}$, see e.g. \cite{BT}, \cite{DKT}, \cite{KKS}, \cite{Ni1}, \cite{Ni2}, 
\cite{NTT} and the references there. For fixed $r\ge3$, we refer to \cite[Theorem 7]{BV}. \\

More generally, let $\gamma_{1}$ and $\gamma_{2}$ be non-zero algebraic numbers satisfying 
$|\gamma_{1}|>1, \gamma_{1}\gamma_{2}=\delta=\pm1$, and define
\[
R_{n}:=g_{1}\gamma_{1}^{n}+g_{2}\gamma_{2}^{n},\ S_{n}:=h_{1}\gamma_{1}^{n}+h_{2}\gamma_{2}^{n}
\quad (n\in\mathbb{Z}),
\]
where $g_{1},g_{2},h_{1},h_{2}$ are non-zero algebraic numbers. Moreover, let 
$\Omega:=(g_{1}h_{2})/(g_{2}h_{1})$. In particular, the special choice 
$(R_{n})=(F_{n}), (S_{n})=(L_{n})$ yields $\Omega=-1$. Next, for fixed $k\in\mathbb{N}$, we introduce 
for $\ell\in\mathbb{Z}, r\in\mathbb{N}\setminus\!\{1\}, \mu=1,\ldots,r-1$ the three series
\begin{equation}\label{a4}
Q_{\mu,r}:=\sum_{h=0}^{\infty}\frac{a_{h}}{\gamma_{1}^{\mu kr^{h}}},\ R_{\ell,r}:=
{\sum_{h=0}^\infty}\,\strut'\frac{b_{h}}{R_{kr^{h}+\ell}},\ S_{\ell,r}:={\sum_{h=0}^\infty}\,
\strut'\frac{c_{h}}{S_{kr^{h}+\ell}},
\end{equation}
where $(a_{h})_{h\ge0}, (b_{h})_{h\ge0}$, and $(c_{h})_{h\ge0}$ are, as above, non-zero periodic 
sequences of algebraic numbers. Then we have the following generalization of Theorem 1.\\

\textbf{Theorem 2}. \textit{If $\Omega\notin(\gamma_{1}/\gamma_{2})^{\mathbb{Z}}$, then the numbers
\begin{equation}\label{a5}
Q_{\mu,r},\ R_{\ell,r},\ S_{\ell,r} \quad (\ell\in\mathbb{Z}, r\in\mathbb{N}\setminus\!\{1\}, 
\mu=1,\ldots,d(r)-1)
\end{equation}
are algebraically independent, unless there exists an $\ell_0\in\mathbb{Z}$ satisfying one of the
following properties: \\
1) $R_{\ell_0}=1$ and $(b_h)$ is constant; \\
2) $S_{\ell_0}=1$ and $(c_h)$ is constant; \\ 
3) $(\gamma_1/\gamma_2)^{\ell_0}=-g_1^{-1}g_2\theta$ and $\Omega\theta\in 
(\gamma_1/\gamma_2)^\mathbb{Z}$, where $\theta=e^{2\pi i/3}$; \\  
4) $(\gamma_1/\gamma_2)^{\ell_0}=-g_1^{-1}g_2\theta$ and $\Omega\theta\in 
(\gamma_1/\gamma_2)^\mathbb{Z}$, where $\theta=e^{-2\pi i/3}$. \\    
In case 1) (or 2)), the number $R_{\ell_0,2}$ (or $S_{\ell_0,2}$, respectively) is algebraic 
but all other numbers in (\ref{a5}) are algebraically independent. In cases 3) and 4), after removing all 
numbers $R_{\ell_0,2^j} \ (j\in\mathbb{N})$ from (\ref{a5}), the remaining numbers are algebraically 
independent.} \\

Again here there are several results considering separately the numbers $R_{\ell,r}$ or $S_{\ell,r}$, 
and even the powers of these, see \cite{NTT}. We would like to point out that also our considerations 
below could be generalized to prove the algebraic independence of
\[
Q_{\mu,r},\ R_{\ell,r,m}:={\sum_{h=0}^\infty}\,\strut'\frac{b_{h}}{(R_{kr^{h}+\ell})^{m}},\ S_{\ell,r,m}:=
{\sum_{h=0}^\infty}\,\strut'\frac{c_{h}}{(S_{kr^{h}+\ell})^{m}}
\]
\[
(\ell\in \mathbb{Z}, m\in \mathbb{N}, r\in \mathbb{N}\setminus\!\{1\}, r\not\in 2^{2\mathbb{N}-1}, 
\mu=1, \ldots, d(r)-1)\ .
\]
The extra condition $r\not\in2^{2\mathbb{N}-1}$ is needed, since in the case $r=2$ there are 
dependence relations with even values of $m$, see \cite[Theorem 2 and Remark 5]{BV1}. \\

Let now $\alpha_\ell \ (\ell\in\mathbb{Z})$ be non-zero algebraic numbers satisfying $|\alpha_{0}|=1, 
|\alpha_{\ell_{1}}|\ne|\alpha_{\ell_{2}}|$ for $\ell_{1}\ne\ell_{2}$, and let $\beta_{\ell}=
\delta_{\ell}\alpha_{\ell}$ with algebraic $\delta_\ell$ such that $|\delta_{\ell}|=1, \delta_{\ell}\ne1
\ (\ell\in \mathbb{Z})$. In proving Theorems 1 and 2, we obtain also a result on the values of the functions
\[
\Gamma_{\mu,r}(z):=\sum_{h=0}^{\infty}a_{h}z^{\mu r^{h}}\quad (\mu=1,\ldots,r-1),
\]
\[
\Phi_{\ell,r}(z):=\sum_{h=0}^{\infty}b_{h}\frac{z^{r^{h}}}{z^{2r^{h}}-\alpha_{\ell}}, \; 
\Lambda_{\ell,r}(z):=\sum_{h=0}^{\infty}c_{h}\frac{z^{r^{h}}}{z^{2r^{h}}-\beta_{\ell}}\quad 
(\ell\in\mathbb{Z}, r\in\mathbb{N}\setminus\!\{1\}).
\] 
\vskip0,3cm

\textbf{Theorem 3}. \textit{Let $\alpha_\ell, \beta_\ell \ (\ell \in \mathbb{Z})$ be as above. If 
$\alpha$ is an algebraic number with $0<\left|\alpha\right|<1$ and $\alpha^{2r^h}\ne\alpha_\ell, 
\beta_\ell$ for all $h\in\mathbb{N}_0$ and $\ell\in\mathbb{Z}$, then the numbers
\begin{equation}\label{a6}
\Gamma_{\mu,r}(\alpha),\ \Phi_{\ell,r}(\alpha),\ \Lambda_{\ell,r}(\alpha) \quad (\ell\in\mathbb{Z},\ 
r\in\mathbb{N}\setminus\!\{1\},\ \mu=1,\ldots,d(r)-1)
\end{equation}
are algebraically independent, unless one of the following cases holds: \\
1) $(b_h)$ is constant and $\alpha_0=1$; \\
2) $(c_h)$ is constant and $\beta_0=1$; \\
3) $\alpha_0=e^{2\pi i/3}$ and $\beta_0=\alpha_0^2$; \\
4) $\alpha_0 = e^{-2\pi i/3}$ and $\beta_0=\alpha_0^2$. \\
In case 1) (or 2)), $\Phi_{0,2}(\alpha)$ (or $\Lambda_{0,2}(\alpha)$, respectively) is algebraic, 
but all other numbers in (\ref{a6}) are algebraically independent. In case 3) (or 4)), after removing 
the numbers $\Phi_{0,2^{j}}(\alpha)$ (or $\Lambda_{0,2^{j}}(\alpha)$) with $j\in\mathbb{N}$ from 
(\ref{a6}) the remaining numbers are algebraically independent.}

\textit{The preceding statement also applies assuming $|\alpha_{0}|=1, 
|\alpha_{\ell_{1}}|\ne|\alpha_{\ell_{2}}|$ and $|\beta_{\ell_{1}}|\ne|\beta_{\ell_{2}}|$ for  
$\ell_{1}\ne\ell_{2}$, and $|\alpha_{\ell_{1}}|\ne|\beta_{\ell_{2}}|$ for all $\ell_1, \ell_2$. Here the 
cases 3) and 4) cannot occur.} \\

Finally, we note the following corollary of Theorem 3 considering the values of Fredholm-type series. \\

\textbf{Corollary 1.} \textit{For algebraic $\alpha$ with $0<|\alpha|<1$, the numbers
\[
\Gamma_{\mu,r}(\alpha)\quad (r\in\mathbb{N}\setminus\!\{1\}, \mu=1,\ldots,d(r)-1)
\]
are algebraically independent}.\\

This gives some further information on the values of series of this type, for earlier results see 
\cite[Chapter 3]{Ni0} and \cite[Theorem 2]{Ni2}. \\

The main ideas of the proofs of our results are similar to \cite{NTT}. First, in Section 2, we consider 
linear independence of Mahler-type functions related to our results. Here the considerations needed 
for the case $\Omega=\Delta(\gamma_1/\gamma_2)^{\ell_1}$ with some $\ell_1\in\mathbb{Z}$ and 
$|\Delta|=1, \Delta\ne1$, are most interesting and challenging, in fact the case 
$\left|\Omega\right|\notin\left|\gamma_1/\gamma_2\right|^{\mathbb{Z}}$ follows essentially from the 
studies in \cite{NTT}. Then, following the lines of \cite{NTT}, the algebraic independence of the functions 
(over $\mathbb{C}(z)$) under consideration and of their values (over $\mathbb{Q}$) are proved, giving 
Theorem 3. In the final Section 5, the connection of the numbers in our theorems with such function 
values implies the validity of Theorems 1 and 2.

\section{Linear independence of functions}

In our first lemma, we are interested in the Mahler-type functions
\[
F_{i}(z):=\sum_{h=0}^{\infty}a^{h}\frac{A_{i}(z^{r^{h}})}{B_{i}(z^{r^{h}})} \quad (i=0,1, \ldots, m),
\]
where $r\geq 2$ is an integer, $a\neq 0$ is a complex number, $A_{i}(z), B_{i}(z)\in \mathbb{C}[z]
\setminus\!\{0\}, A_{i}(0)= 0, B_{0}(z)\equiv 1$, and, for any $i\geq 1, A_{i}(z)$ and $B_{i}(z)$ are 
coprime, and the $B_{i}(z)$ are distinct, non-constant, and monic. Clearly
\[
aF_{i}(z^{r})=F_{i}(z)-\frac{A_{i}(z)}{B_{i}(z)} \quad (i=0,1, \ldots, m).
\]

\textbf{Lemma 1}. \textit{For each $i=1,\ldots,m$, assume that $B_i(z)$ (as above) has exactly $t_i$ 
distinct zeros all of them of the same absolute value $\omega_i$. Let all $\omega_1,\ldots,\omega_m$
be distinct, and $\omega_m=1$. In case $m\ge2$, suppose furthermore $(r-1)t_i\ge t_j$ for every 
pair $(i,j)$ with $1\le i,j<m$. Then the function 
\[
u_0F_0(z)+u_1F_1(z)+\cdots+u_mF_m(z) 
\]
with $(u_0,\ldots,u_m)\in\mathbb{C}^{m+1}$ is rational if and only if $u_1=\cdots=u_{m-1}=0$ and 
the function $u_0F_0(z)+u_mF_m(z)$ is rational.} \\

\textit{Proof}. Clearly, we may assume $m\ge2$ and $(u_1,\ldots,u_{m-1})\ne\underline{0}$. The
function $g(z):=\sum_{i=0}^m u_iF_i(z)$ in our lemma satisfies the functional equation
\begin{equation}\label{a7}
a g(z^r)=g(z)-u_0 A_0(z)-u_1\frac{A_1(z)}{B_1(z)}-\cdots-u_m\frac{A_m(z)}{B_m(z)}.
\end{equation}
Assuming that $g(z)$ is rational, it must have poles of absolute value $\ne1$, and we may first 
suppose of absolute value 
$>1$. Let $p$ be such a pole with maximal absolute value. Then $|p|=\omega_i$ with some 
$i\in\{1,\ldots,m-1\}$. Clearly, $p$ cannot be a pole of $g(z^r)$. Thus, the hypotheses of our lemma 
together with (\ref{a7}) imply that all distinct zeros of $B_i(z)$, say $\alpha_{i,1},\ldots,
\alpha_{i,t_i}$, are poles of $g(z)$, hence $r t_i$ numbers $\sqrt[r]{\alpha_{i,\nu}}$ are poles of 
$g(z^r)$. By the distinctness of the $\omega$'s, the function $ag(z^r)-g(z)$ has exactly $t_j$ 
poles of the same absolute value $\omega_j$ (if $u_j\ne0$), and therefore at least $r t_i-t_j 
\ (\ge t_i)$ of the above $\sqrt[r]{\alpha_{i,\nu}}$ are poles of $g(z)$ (assuming $\sqrt[r]{\omega_i}
=\omega_j$). Let these be $q_1,\ldots,q_v$ with $v\ge t_i$. The $r v$ numbers $\sqrt[r]{q_i}$ 
are poles of $g(z^r)$, and again at least $r v-t_k \ (\ge t_i)$ of these are poles of $g(z)$. By 
repeating this conclusion, we get a contradiction. The same argument works if $g(z)$ has poles of 
absolute value $<1$ (but note that 0 is not a pole of $g(z)$). 
\hfill $\Box$ \\

We now consider some special functions of the above type, namely
\begin{equation}\label{a8}
\gamma_{\mu}(a,z):=\sum_{h=0}^{\infty}a^{h}z^{\mu r^{h}}\quad (\mu=1, \ldots, r-1),
\end{equation}
\vspace*{-0.5cm}
\begin{equation}\label{a9}
\varphi_{\ell}(a,z):=\sum_{h=0}^{\infty}a^{h}\frac{z^{r^{h}}}{z^{2r^{h}}-\alpha_{\ell}},\enspace
\lambda_{\ell}(a,z)=\sum_{h=0}^{\infty}a^{h}\frac{z^{r^{h}}}{z^{2r^{h}}-\beta_{\ell}}\quad (\ell\in\mathbb{Z}),
\end{equation}
where $a$ and all $\alpha_\ell,\beta_{\ell}$ are non-zero complex numbers. These functions 
satisfy the following functional equations
\[
a\gamma_{\mu}(a,z^{r})=\gamma_{\mu}(a,z)-z^{\mu},\enspace a\varphi_{\ell}(a,z^{r})
=\varphi_{\ell}(a,z)-\frac{z}{z^{2}-\alpha_{\ell}},
\]
\vspace*{-1.0cm}
\begin{equation}\label{a10}
{}
\end{equation}
\vspace*{-1.0cm}
\[
a\lambda_{\ell}(a,z^{r})=\lambda_{\ell}(a,z)-\frac{z}{z^{2}-\beta_{\ell}}.
\]

\textbf{Lemma 2}. \textit{Suppose $|\alpha_{0}|=1, |\alpha_{\ell_{1}}|\ne|\alpha_{\ell_{2}}|$ for 
$\ell_{1}\ne\ell_{2}$, and $\beta_{\ell}=\delta_{\ell}\alpha_{\ell}$ with $|\delta_{\ell}|=1, \delta_{\ell}\ne1$
for any $\ell\in\mathbb{Z}$. If $p_1,\ldots,p_{r-1}$ and, for some non-negative integer $L$, the
$u_{\ell},v_{\ell}$ with $|\ell|\le L$ are complex numbers, then the function
\begin{equation}\label{a11}
\sum_{\mu=1}^{r-1}p_{\mu}\gamma_{\mu}(a,z)+\sum_{\ell=-L}^{L}\Big(u_\ell\varphi_{\ell}(a,z)+
v_{\ell}\lambda_{\ell}(a,z)\Big)
\end{equation}
is rational if and only if $u_{\ell}=v_{\ell}=0$ for any $\ell\ne0$, and the following function is rational}
\begin{equation}\label{a12}
\sum_{\mu=1}^{r-1}p_{\mu}\gamma_{\mu}(a,z)+u_{0}\varphi_{0}(a,z)+v_{0}\lambda_{0}(a,z).
\end{equation}

\textit{Proof}. If $g(z)$ denotes the function (\ref{a11}), then
\begin{equation}\label{p1}
ag(z^{r})=g(z)-A(z)-\sum_{\ell=-L}^{L}\Big(u_\ell\frac{z}{z^{2}-\alpha_{\ell}}+v_{\ell}\frac{z}{z^{2}-
\beta_{\ell}}\Big)
\end{equation}
holds, by (\ref{a10}), with $A(z):=\sum_{\mu=1}^{r-1}p_\mu z^\mu$. Assume now that $g(z)$ is rational. 

We begin with the case $r\ge3$, where we intend to apply Lemma 1. Consider the integers $\ell$
with $0<|\ell|\le L$. If $u_{\ell}v_{\ell}=0$ holds for such an $\ell$ but $(u_\ell,v_\ell)\ne(0,0)$, then
the corresponding $t_\ell$ in the sense of Lemma 1 equals 2. If, on the other hand, $u_{\ell}v_{\ell}\ne0$,
then the contribution 
\[
u_\ell\frac{z}{z^{2}-\alpha_{\ell}}+v_{\ell}\frac{z}{z^{2}-\beta_{\ell}}= \frac{P_\ell(z)}{(z^{2}-\alpha_{\ell})
(z^{2}-\beta_{\ell})} \enspace \mathrm{with} \enspace P_\ell(z):=z\big((u_\ell+v_\ell)z^2-
(u_\ell\beta_\ell+v_\ell\alpha_\ell)\big)
\]
to the sum in (\ref{p1}) has coprime non-zero $P_\ell(z), (z^{2}-\alpha_{\ell})(z^{2}-\beta_{\ell})$,
whence the corresponding $t_\ell$ is 4. With these two possibilities for the $t$'s in Lemma 1,
we see that the condition $(r-1)t_i\ge t_j$ imposed there is satisfied, by $r\ge3$. Thus, Lemma 1
implies Lemma 2 in this case.

We are left with the case $r=2$. Suppose that there exists some $\ell, \;0<|\ell|\le L,$ with
$(u_\ell,v_\ell)\ne(0,0)$ and, moreover, $|\alpha_{\ell}|>1$. Let $\ell=m$ be the index with smallest 
$|\alpha_{\ell}|$ having this property.

Before arguing further, we need the following intermediate claim. If $\omega$ with $|\omega|>1$ is a 
pole of $g(z)$, then $|\omega|\ge|\alpha_{m}|$. Assume, on the contrary, that $1<|\omega|<|\alpha_{m}|$. 
Then $\omega_{1}:=\sqrt{\omega}$ is a pole of $g(z^{2})$ and also a pole of $g(z)$, by
\begin{equation}\label{b1}
ag(z^{2})=g(z)-p_1z-\sum_{\ell=-L}^{L}\Big(u_\ell\frac{z}{z^{2}-\alpha_{\ell}}+v_{\ell}\frac{z}{z^{2}-
\beta_{\ell}}\Big).
\end{equation}
But then $\omega_{2}:=\sqrt{\omega_{1}}$ is a pole of $g(z^{2})$, hence a pole of $g(z)$, again by 
(\ref{b1}). Repeating this argument we get a contradiction proving our claim.

With $m$ defined above, we assume $u_m\ne0$ (the case $v_m\ne0$ being similar). As we saw,
$g(z)$ is polefree in the circular annulus $1<|z|<|\alpha_m|$. Hence, by (\ref{b1}), the
function
\[
h(z):=ag(z^2)+\frac{u_mz}{z^2-\alpha_m}
\]
is holomorphic in our annulus from which we have removed the points
$\sqrt{\beta_m}$ and $-\sqrt{\beta_m}$ (if $v_m\ne0$), and the points 
$\pm\sqrt{\alpha_\ell}, \pm\sqrt{\beta_\ell} \ (\ell\ne m)$ if these points lie in the annulus. Letting now 
$z\to\sqrt{\alpha_m}$
(or $z\to-\sqrt{\alpha_m})$, one sees from the definition of $h(z)$ that $g(z)$ must have a pole
at $\alpha_m$. Assume that $b(z-\alpha_m)^{-k}+\cdots$ with $b\ne0$ and a positive integer $k$
indicates the beginning of the Laurent expansion of $g(z)$ at $\alpha_m$, we find
\[
h(z)=\Big(\frac{ab}{(z^2-\alpha_m)^k}+\cdots\Big)+\frac{u_mz}{z^2-\alpha_m}
\]
near $\sqrt{\alpha_m}$ and $-\sqrt{\alpha_m}$. Thus, we obtain $k=1$ and $h(z)=(ab+u_mz)/(z^2-\alpha_m)
+\cdots$, hence $ab+u_m\sqrt{\alpha_m}=0$ and $ab-u_m\sqrt{\alpha_m}=0$ leading to $u_m=0$,
a contradiction. Thus, Lemma 2 is proved if there exists some $\ell, 0<|\ell|\le L$, with
$(u_\ell,v_\ell)\ne(0,0)$ and $|\alpha_\ell|>1$. If, for all $\ell, 0<|\ell|\le L$, with $(u_\ell,v_\ell)\ne(0,0)$, 
the inequality $|\alpha_\ell|<1$ holds, the proof is similar. \hfill $\Box$ \\

\textit{Remark 1}. If we assume that $|\alpha_{0}|=1, |\alpha_{\ell_{1}}|\ne|\alpha_{\ell_{2}}|$ and 
$|\beta_{\ell_{1}}|\ne|\beta_{\ell_{2}}|$ for $\ell_{1}\ne\ell_{2}$, and 
$|\alpha_{\ell_{1}}|\ne|\beta_{\ell_{2}}|$ for all $\ell_1, \ell_2$, then we see, as at the beginning of 
the above proof of Lemma 2, that the function $g(z)$ in (\ref{a11}) is rational if and only if 
$v_0=0, u_\ell=v_\ell=0$ for all $\ell\ne0$, and the following function is rational
\begin{equation}\label{b2}
\sum_{\mu=1}^{r-1}p_{\mu}\gamma_{\mu}(a,z)+u_{0}\varphi_{0}(a,z).
\end{equation}

The subsequent Lemma 3 concerns just the question in which cases the function (\ref{a12}) can be
rational.  \\

\textbf{Lemma 3}. \textit{Suppose $|\alpha_0|=1$ and $\beta_0=\delta_0\alpha_0$ with $|\delta_0|=1,
\delta_0\ne1$. Suppose that $p_1,\ldots,p_{r-1},u_0,$\ $v_0$ are complex numbers, not all zero,
and let $g_0(z)$ denote the function (\ref{a12}). If either $r\ge4$, or $r=3$ and $a^2\ne9$, then 
$g_0(z)$ is not rational. If $r=2$ and $a^2\ne2, 4$, then this function is rational 
if and only if $p_1=0$ and, moreover, one of the following six conditions holds, $\zeta$
denoting} $e^{\pi i/3}$: \\
1) $a=1, \alpha_{0}=1, v_{0}=0$; \\
2) $a=1, \beta_{0}=1, u_{0}=0$; \\
3,4) $a=\pm1, \alpha_0=\zeta^2, \beta_0=\zeta^4, v_0=-a u_0\zeta$; \\
5,6) $a=\pm1, \alpha_0=\zeta^4, \beta_0=\zeta^2, u_0=-a v_0\zeta$. \\

\textit{Remark 2.}
In case 1), $\varphi_0(1,z)$ is the rational function $z/(z-1)$, and in case 2), $\lambda_0(1,z)$  is
this same function. In case 3),
\[
\zeta^{-2}\varphi_0(1,z)+\zeta^2\lambda_0(1,z)=\frac{2z^2+z}{z^2+z+1}
\]
holds, and in case 5), the roles of $\varphi_0(1,z)$ and $\lambda_0(1,z)$ are reversed. In case 4),
\[
\frac{1}{1+\zeta}\varphi_0(-1,z)+\frac{\zeta}{1+\zeta}\lambda_0(-1,z)=\frac{z}{z^2+z+1}
\]
holds, whereas in case 6), the roles of $\varphi_0(1,z)$ and $\lambda_0(1,z)$ are reversed. But it should
be noted that, in all cases concerning $a=1$, a solution $g$ of (\ref{p1}), (\ref{b1}), or (\ref{a13})
below is only determined up to an arbitrary additive constant. \\

\textit{Proof of Lemma 3.} Our function $g_0(z)$ satisfies
\begin{equation}\label{a13}
ag_{0}(z^{r})=g_{0}(z)-A(z)-u_{0}\frac{z}{z^{2}-\alpha_{0}}-v_{0}\frac{z}{z^{2}-\beta_{0}},
\end{equation}
where $A(z)$ is the polynomial from the beginning of the proof of Lemma 2. 

We first show: If $(u_0,v_0)=(0,0)$ (implying $A(z)\ne0$), then $g_0(z)$ is not rational.
Assuming, on the contrary, that $g_0(z)=U(z)/V(z)$ with coprime $U(z),V(z)\in\mathbb{C}[z]
\!\setminus\!\{0\}$, we have
\[
aU(z^r)V(z)=\big(U(z)-A(z)V(z)\big)V(z^r)
\]
hence $V(z^r)|V(z)$, and therefore $\deg V=0$, whence $g_0(z)$ is a non-constant polynomial
with $G:=\deg g_0 \;(\ge1)$, say. Hence $\deg(ag_0(z^r)=Gr, \deg(g_0(z)-A(z))\le\max(G,r-1)$
giving the desired contradiction. (Note here: If $g_0$ would be constant, then $A$ as well, by
(\ref{a13}) with $(u_0,v_0)=(0,0)$.)

Assume now that $u_{0}\ne0$ and $v_{0}=0$, the case $v_{0}\ne0, u_{0}=0$ being similar. Clearly,
all poles $\omega$ of $g_{0}(z)$ satisfy $|\omega|=1$, and $rt\le t+2$ or $t\le2/(r-1)$, where $t$ 
denotes the number of distinct poles. Thus, we have a contradiction if $r\ge4$. If $r=3$, then 
\cite[Lemma 2]{BV} gives a contradiction. In case $r=2$, \cite[Lemma 3]{BV} implies that the only 
possibility is $p_1=0, a=\alpha_0=1$ and, in this case, $\varphi_{0}(1,z)$ is the rational function
we already saw in Remark 2. \\

Next, let us assume that $u_{0}v_{0}\ne0$. This time we get the upper bound $t\le4/(r-1)$ for the 
number of distinct poles of $g_0(z)$. Thus, we need to consider only the values $r=2,3,4$ and 5,
and $t\le4$ if $r=2, t\le2$ if $r=3$, and $t=1$ otherwise. By comparing the arguments of the possible 
poles on both sides of (\ref{a13}), we immediately get a contradiction if $r=5$. In case $r=4$,
\cite[Lemma 3]{BV} gives a contradiction, too. \\

Assume now $r=3$. Then $t\le2$, and we immediately see that we must have $t=2$. Let 
$\omega_{1}$ and $\omega_{2}$ be the distinct poles of $g_0(z)$ with arguments $0\le\tau_1<\tau_2
<2\pi$, respectively. The arguments of the poles on the left-hand side of (\ref{a13}) are 
$\tau_1/3+j2\pi/3, \tau_2/3+j2\pi/3, j=0,1,2$, and the arguments of the possible poles on the right-hand 
side are $\tau_1, \tau_2, \phi_1/2, \phi_1/2+\pi, \phi_2/2, \phi_2/2+\pi$, where $\phi_1=
\min(\arg\alpha_0, \arg\beta_0), \phi_2=\max(\arg\alpha_0, \arg\beta_0),  0\le\phi_1<\phi_2<2\pi$. 
If $\tau_1=0$, then comparision of the arguments of the poles gives us one possibility, where 
$\phi_1=2\pi/3, \phi_2=4\pi/3, \tau_1=0, \tau_2=\pi$. By denoting $\zeta=e^{\pi i/3}$, we get $z^3-1
=(z-1)(z-\zeta^2)(z-\zeta^4), z^3+1=(z+1)(z-\zeta)(z-\zeta^5), z^2-\alpha_0=(z-\zeta)(z-\zeta^4), 
z^2-\beta_0=(z-\zeta^2)(z-\zeta^5)$ , where we assumed (w.l.o.g.) $\phi_1=\arg\alpha_0, \phi_2
=\arg\beta_0$. Therefore $g_0(z)$ in (\ref{a13}) takes the shape $g_0(z)=P(z)/(z^2-1)$, where
$P(\pm1)\ne0$, and we get
\[
aP(z^3)=P(z)(z^4+z^2+1)-A(z)(z^6-1)-u_0z(z^2-1)(z^2-\zeta^4)-v_0z(z^2-1)(z^2-\zeta^2).
\]
This yields $aP(1)=3P(1)$, and since $P(1)\ne0$, this contradicts our assumption $a\ne3$.

In the case $\tau_1>0$, we have also one possibility, $\phi_1=\pi/3, \phi_2=5\pi/3, \tau_1=\pi/2, 
\tau_2=3\pi/2$. Thus $g_0(z)=P(z)/(z^2+1)$, where $P(\pm i)\ne0$. By (\ref{a13}),
\[
aP(z^3)=P(z)(z^4-z^2+1)-A(z)(z^6+1)-u_0z(z^2+1)(z^2-\zeta^5)-v_0z(z^2+1)(z^2-\zeta).
\]
Thus, $a^2P(i)=3aP(-i)=9P(i)$ giving $a^2=9$, which contradicts our assumptions. \\ 

We are left with the case $r=2,\, t\le4$ as noticed above. We assume firstly that none of
$\alpha_0,\beta_0$ is equal to 1. \\

Let $t=4$ and $\omega_j$ with $0\le\tau_j=\arg\omega_j<2\pi \ (j=1,2,3,4)$ be the poles. W.l.o.g. we 
may assume that $0\le\tau_1<\tau_2<\tau_3<\tau_4<2\pi$ and $0<\phi_1=\arg\alpha_0<\phi_2=
\arg\beta_0<2\pi$. Comparing the arguments on both sides of (\ref{a13}) gives $\tau_4=\phi_2$. 
If $\tau_1=0$, then we must have $\tau_2=\phi_1$ and $\{\tau_3/2, \pi, \tau_3/2+\pi\}=
\{\phi_1,\ \tau_3,\ \phi_2\}$. This means that $\tau_3=\pi, \phi_1=\pi/2, \phi_2=3\pi/2$, and $\omega_1=1, 
\omega_2=\alpha_0=i, \omega_3=-1, \omega_4=\beta_0=-i$. Thus, $g_0(z)=P(z)/(z^4-1)$, where 
$P(z)$ is a polynomial satisfying $P(\pm1)P(\pm i)\ne0$. Now (\ref{a13}) gives
\[
aP(z^2)=P(z)(z^4+1)-p_1z(z^4-1)(z^4+1)-u_0z(z^4-1)(z^2+i)-v_0z(z^4-1)(z^2-i),
\]
hence $aP(1)=2P(1)$, and so $a=2$, contrary to our assumptions. In the case $\tau_1>0$, we must 
have $\tau_1=\phi_1$ and $\{\tau_2/2, \tau_3/2, \tau_2/2+\pi, \tau_3/2+\pi\}=\{\phi_1, \tau_2, \tau_3, 
\phi_2\}$. We therefore have $\tau_2=2\phi_1, \tau_3=4\phi_1, \phi_1+\pi=4\phi_1$ or $\phi_1=\pi/3$, 
and finally $\phi_2=5\pi/3$. By denoting $\zeta=e^{\pi i/3}$ as above, we have $\omega_1=\alpha_0
=\zeta, \omega_2=\zeta^2, \omega_3=\zeta^4, \omega_4=\beta_0=\zeta^5$. As before, we now have
\[
aP(z^2)=P(z)(z^2-\zeta)(z^2+\zeta^2)-p_1z(z^4-\zeta^2)(z^4-\zeta^4)-u_0z(z^2+\zeta)(z^4-\zeta^4)
-v_0z(z^4-\zeta^2)(z^2-\zeta^2),
\]
where $P(z)$ satisfies $P(\pm\zeta)P(\pm\zeta^2)\ne0$. On substituting $z=-\zeta$ and $z=\zeta^2$
to the last equation, we are led to
\[
aP(\zeta^2)=-2\zeta^2P(-\zeta) \enspace \mathrm{and} \enspace aP(-\zeta)=2\zeta P(\zeta^2),
\]
respectively, implying $a^2=4$ but again this contradicts our assumptions. \\

Let now $t=3$. Comparing the arguments of the poles gives here four possibilities: 1) $\omega_1=1, 
\omega_2=\alpha_0=i, \omega_3=\beta_0=-1$; 2) $\omega_1=1, \omega_2=\alpha_0=-1, \omega_3
=\beta_0=-i$; 3) $\omega_1=\alpha_0=\zeta, \omega_2=\zeta^2, \omega_3=\beta_0=\zeta^4$; 
4) $\omega_1=\alpha_0=\zeta^2, \omega_2=\zeta^4, \omega_3=\beta_0=\zeta^5$.

In case 1), we apply (\ref{a13}) and obtain
\[
aP(z^2)=P(z)(z+i)(z^2-i)-p_1z(z^4-1)(z^2-i)-u_0z(z^4-1)-v_0z(z^2-1)(z^2-i)
\]
with some polynomial $P(z)$ satisfying $P(\pm 1)P(i)\ne0$. Therefore $aP(1)=2P(1)$ and $a=2$. 
In case 2), we get the same contradiction to our assumptions.

Using (\ref{a13}) we get in case 3)
\[
aP(z^2)=P(z)(z^2-\zeta)(z+\zeta^2)-p_1z(z^4-\zeta^2)(z^2-\zeta^2)-u_0z(z^2-\zeta^2)(z^2+\zeta)
-v_0z(z^2-\zeta)(z^2-\zeta^2),
\]
where $P(z)$ satisfies $P(\pm\zeta)P(\zeta^2)\ne0$. Taking here $z=\zeta^2$ and $z=-\zeta^2$ we obtain
\[
aP(-\zeta)=4P(\zeta^2)+2v_0\zeta(\zeta+1) \enspace \mathrm{and} \enspace aP(-\zeta)=-2v_0\zeta(\zeta+1),
\]
respectively, hence $aP(-\zeta)=2P(\zeta^2)$. Inserting, moreover, $z=-\zeta$ we find
$aP(\zeta^2)=P(-\zeta)$, and the last two equations lead to $a^{2}=2$, contrary to our assumptions. 
Case 4) can be similarly dealt with. \\

Assume finally $t=2$. In this case, a reasoning as above implies one possibility, namely $\omega_1
=\alpha_0=\zeta^2, \omega_2=\beta_0=\zeta^4$. By using (\ref{a13}), we therefore get
\begin{equation}\label{b3}
aP(z^2)=P(z)(z-\zeta)(z+\zeta^2)-p_1z(z^2-\zeta^2)(z^2-\zeta^4)-u_0z(z^2-\zeta^4)-v_0z(z^2-\zeta^2),
\end{equation}
where $P(z)$ satisfies $P(\zeta^2)P(\zeta^4)\ne0$. For $z=\zeta$ and $z=-\zeta^2$, we get $aP(\zeta^2)=
-u_0\zeta^2(1+\zeta)$ and $aP(\zeta^4)=v_0(1+\zeta)$, respectively. Thus, applying (\ref{b3}) with
$z=\zeta^2$ and $z=-\zeta$, we find after some minor computation
\[
aP(\zeta^4)=-\zeta^2P(\zeta^2) \enspace \mathrm{and} \enspace aP(\zeta^2)=\zeta P(\zeta^4),
\]
respectively. This gives immediately $a^2=1$ hence $a=\pm1$.

If $a=-1$, then the above implies $u_0\zeta^2(1+\zeta)=P(\zeta^2)=-\zeta P(\zeta^4)
=v_0\zeta(1+\zeta)$. Thus, $v_0=u_0\zeta$ and (\ref{b3}) gives
\[
-P(z^2)=P(z)(z^2-z+1)-p_1z(z^2-\zeta^2)(z^2-\zeta^4)-u_0(z(z^2-\zeta^4)+\zeta z(z^2-\zeta^2)).
\]
From this equation, we deduce $\deg P(z)\le2$, say $P(z)=P_1z+P_2z^2$, and $p_1=0$. Then the 
above equation leads to $P_1=u_0(1+\zeta), P_2=0$, and the rational function
\[
g_0(z)=u_0(1+\zeta)\frac{z}{z^2+z+1}=u_0\varphi_0(-1,z)+u_0\zeta\lambda_0(-1,z)
\]
satisfies (\ref{a13}) in this case. 

In case $a=1$, we have $-u_0\zeta^2(1+\zeta)=P(\zeta^2)=\zeta P(\zeta^4)=v_0\zeta(1+\zeta)$ 
giving $v_0=-u_0\zeta$. By (\ref{b3}),
\[
P(z^2)=P(z)(z^2-z+1))-p_1z(z^2-\zeta^2)(z^2-\zeta^4)-u_0(z(z^2-\zeta^4)-\zeta z(z^2-\zeta^2)).
\]
Again we are led to $\deg P(z)\le2$ and $p_1=0$. By substituting $P(z)=P_1z+P_2z^2$ 
to the above equation, we obtain $P_1=-u_0(1-\zeta), P_2=2P_1$. (Note here that, in virtue of our 
concluding statement in Remark 2, we are allowed to assume $P(0)=0$, w.l.o.g.) Thus, in this case, 
(\ref{a13}) has a rational solution
\[
g_0(z)=u_0(\zeta-1)\frac{2z^2+z}{z^2+z+1}=u_0\varphi_0(1,z)-u_0\zeta\lambda_0(1,z).
\]     

We still have to consider the case $\alpha_0=1$. Now the cases $t=3$ or $t=4$ are clearly not 
possible, and $t=2$ is only possible if $\omega_1=1, \omega_2=\beta_0=-1$. By (\ref{a13}),
\[
aP(z^2)=P(z)(z^2+1)-p_1z(z^2-1)(z^2+1)-u_0z(z^2+1)-v_0z(z^2-1)
\]
with some polynomial $P(z)$ satisfying $P(\pm1)\ne0$. Therefore $aP(-1)=2iv_0$ and $aP(-1)=-2iv_0$ 
giving $P(-1)=0$, a contradiction.

Thus, the condition $u_0v_0\ne0$ holds exactly in cases 3)-6), and this completes the proof of 
Lemma 3. \hfill $\Box$ \\

\textit{Remark 3.} From the beginning of the above proof it follows that, under the assumptions of 
Remark 1, the function (\ref{b2}) is rational if and only if $r=2, p_1=0, a=\alpha_0=1$ and, in this 
case, $\varphi_0(1,z)$ is rational.

\section{Algebraic independence of functions}

For a sequence $(a_h)_{h\ge0}$ of complex numbers, we denote by $(a_h^{(j)})_{h\ge0}, j\in\mathbb{N}$,
the sequences
\[
(a_h^{(1)})=(a_0,a_1,\ldots), (a_h^{(2)})=(a_0,0,a_1,0,a_2,0,\ldots), (a_h^{(3)})=(a_0,0,0,a_1,0,0,a_2,0,0\ldots), 
\ldots.
\]
By \cite[Lemmas 2.7 and 2.8]{NTT}, these sequences have the following properties: \\
a) \textit{If $(a_h)_{h\ge0}$ is a periodic sequence, not identically zero, then the sequences
$(a_h^{(j)})_{h\ge0} \ (j\ge1)$ \\ 
\hspace*{0.3cm} are linearly independent over} $\mathbb{C}$.\\
b) \textit{If $(a_h)_{h\ge0}$ is a periodic sequence with period length greater than 1, then  
$(1),(a_h^{(j)})_{h\ge0} \ (j\ge1)$
\hspace*{0.4cm}are linearly independent over $\mathbb{C}$, $(1)$ being the obvious constant sequence.} \\

The following lemma will be needed (see \cite[p.106]{NTT}). \\

\textbf{Lemma 4}. \textit{Let $p$ be a positive integer,
$(B_h)_{h\ge0}$ a periodic sequence with period lenght dividing $p$, 
and let $R(z)$ be the quotient of two polynomials in $z=(z_1,\ldots, z_n)$ such that the numerator 
vanishes at the origin of $\mathbb{C}^n$ but the denominator does not. Further, define
\[
f_r(z)=\sum_{h=0}^\infty B_hR(z^{r^h}).
\]
Then, for any $s\in\mathbb{N}$,
\[
f_{r^j}(z^{r^{s!p}})=f_{r^j}(z)+R_{r,j}(z)
\]
holds for $j=1,\ldots,s$  with rational functions $R_{r,j}(z).\; ($Note here $z^k=(z_1^k,\ldots,z_n^k))$}. \\

\textit{Proof}. According to the definition of $f_r(z)$, we have $f_{r^j}(z)=\sum_{h\ge0}B_hR(z^{r^{jh}})$, 
hence
\[
f_{r^j}(z^{r^{s!p}})=\sum_{h=0}^\infty B_hR(z^{r^{j(h+(s!/j)p)}})=\sum_{\widetilde{h}=(s!/j)p}^\infty
B_{\widetilde{h}-(s!/j)p}R(z^{r^{j\widetilde{h}}})=\sum_{\widetilde{h}=\widetilde{h}(j)}^\infty 
B_{\widetilde{h}}R(z^{r^{j\widetilde{h}}})
\]
(with $\widetilde{h}(j):=(s!/j)p\in\mathbb{N}$), and this equals $f_{r^j}(z)$ up to a rational function. 
\hfill $\Box$ \\ 

Let now $t_1,t_2,t_3\in\mathbb{N}$, and suppose
\[
q_{k,m}(z)\in \mathbb{C}(z_{1}, \ldots, z_{n}) \quad (k=1,2,3;\ m=1,\ldots,t_k),
\]
where the numerators vanish at the origin but the denominators do not. Then we define
\begin{equation}\label{a14}
f_{k,m}(a,z):=\sum_{h=0}^\infty a^hq_{k,m}(z^{r^h}) \quad (k=1,2,3;\ m=1,\ldots,t_k).
\end{equation}
Further, let $(b_{k,h})_{h\ge0}\; (k=1,2,3)$ be non-zero periodic sequences of complex numbers with 
period lengths $p_k$, respectively, put $p:=\mathrm{lcm}(p_1,p_2,p_3)$, and denote, moreover,
\begin{equation}\label{a15}
f_{k,m,r}(z):=\sum_{h=0}^\infty b_{k,h}q_{k,m}(z^{r^h}) \quad (k=1,2,3;\ m=1,\ldots, t_k).
\end{equation}

By following ideas of \cite[pp.106-107]{NTT}, we can now prove \\

\textbf{Lemma 5}. \textit{If, for any root of unity $\xi$, the functions (\ref{a14}) formed with $a=\xi$ are 
linearly independent over $\mathbb{C}$ modulo $\mathbb{C}(z_1,\ldots,z_n)$ , then the functions
\[
f_{k,m,r^j}(z) \quad (j\in\mathbb{N};\ k=1,2,3;\ m=1,\ldots,t_k)
\]
are algebraically independent over $\mathbb{C}(z_1,\ldots,z_n)$}. \\

\textit{Proof}. Assume, on the contrary, that there exists an $s\in\mathbb{N}$ such that the functions
\[
f_{k,m,r^j}(z) \quad (j=1,\ldots,s;\ k=1,2,3;\ m=1,\ldots,t_k)
\]
are algebraically dependent. By Lemma 4, we may apply \cite[Corollary of Theorem 3.2.1]{Ni0}, 
whence these functions are linearly dependent over $\mathbb{C}$ modulo $\mathbb{C}(z_1,\ldots,z_n)$.
Thus, there exist complex numbers $a_{k,m,j}$, not all zero, such that
\[
G(z):=\sum_{k=1}^3\sum_{m=1}^{t_k}\sum_{j=1}^{s}a_{k,m,j}f_{k,m,r^j}(z)\in\mathbb{C}(z_1,\ldots,z_n).
\]
Since
\[
f_{k,m,r^j}(z)=\sum_{h=0}^\infty b_{k,h}q_{k,m}(z^{r^{jh}})=\sum_{h=0}^\infty b_{k,h}^{(j)}q_{k,m}(z^{r^h}),
\]
we obtain
\[
G(z)=\sum_{k=1}^3\sum_{m=1}^{t_k}\sum_{h=0}^\infty\Big(\sum_{j=1}^sa_{k,m,j}b_{k,h}^{(j)}\Big)
q_{k,m}(z^{r^h}).
\]
Here all sequences
\[
\Big(\sum_{j=1}^sa_{k,m,j}b_{k,h}^{(j)}\Big)_{h\ge0}
\]
are periodic with period lengths dividing $s!p$, and therefore there exist complex numbers $A_{k,m,i}$ 
such that, for all $h\ge0$,
\begin{equation}\label{a16}
\sum_{j=1}^sa_{k,m,j}b_{k,h}^{(j)}=\sum_{i=0}^{s!p-1}A_{k,m,i}\omega^{ih}
\end{equation}
with a primitive $(s!p)$-th root of unity $\omega$. By property a) from the beginning of this section, not 
all $A_{k,m,i}$ vanish. Thus, $G(z)$ has the form
\begin{equation}\label{a17}
G(z)=\sum_{k=1}^3\sum_{m=1}^{t_k}\sum_{i=0}^{s!p-1}A_{k,m,i}f_{k,m}(\omega^i,z)=\sum_{i=0}^{s!p-1}
G_i(z)\in\mathbb{C}(z_1,\ldots,z_n),
\end{equation}
where
\[
G_i(z):=\sum_{k=1}^3\sum_{m=1}^{t_k}A_{k,m,i}f_{k,m}(\omega^i,z).
\]
Here
\[
\omega^iG_i(z^r)-G_i(z)\in\mathbb{C}(z_1,\ldots, z_n).
\]
Let now $J := \{i: 0 \leq i \leq s!p - 1$ and at least one $A_{k,m,i} \neq 0\}$. As noted above $J$ is not empty. By applying again \cite[Corollary of Theorem 3.2.1]{Ni0} to the functions $G_i(z), i \in J$, it follows that there exists an $i_0 \in J$ such that 
$G_{i_0}(z)\in\mathbb{C}(z_1,\ldots,z_n)$ (note that the $\omega^i$ are distinct). 
By our hypothesis on linear independence of the $f_{k,m}(\omega^{i_0},z)$, we get the contradiction that 
all $A_{k,m,i_0}$ vanish, proving Lemma 5. \hfill $\Box$ \\

We now apply the above results to the functions of one variable
\begin{equation}\label{a18}
\gamma_{\mu,r}(z):=\sum_{h=0}^\infty a_hz^{\mu r^h}  \quad (\mu=1,\ldots,r-1),
\end{equation}
\begin{equation}\label{a19}
\varphi_{\ell,r}(z):=\sum_{h=0}^\infty b_h\frac{z^{r^h}}{z^{2r^h}-\alpha_{\ell,r}} \quad 
(\ell\in\mathbb{Z}, r\in\mathbb{N}, r\ge2),
\end{equation}
\begin{equation}\label{a20}
\lambda_{\ell,r}(z):=\sum_{h=0}^\infty c_h\frac{z^{r^h}}{z^{2r^h}-\beta_{\ell,r}} \quad 
(\ell\in\mathbb{Z}, r\in\mathbb{N}, r\ge2).
\end{equation} 
\vskip0.2cm
\textbf{Lemma 6}. \textit{Assume that $\alpha_{\ell,r}$ and $\beta_{\ell,r}$ satisfy the conditions of 
$\alpha_\ell$ and $\beta_\ell$ in Lemma 2, and let $\alpha_{\ell,r^j}:=\alpha_{\ell,r}$ and 
$\beta_{\ell,r^j}:=\beta_{\ell,r}$ for all $j\in\mathbb{N}$. If $r\ge3$, then the functions
\begin{equation}\label{a21}
\gamma_{\mu,r^j}(z),\ \varphi_{\ell,r^j}(z),\ \lambda_{\ell,r^j}(z)  \quad (\mu=1,\ldots,r-1,\
\ell\in\mathbb{Z},\ j\in\mathbb{N})
\end{equation}
are algebraically independent over $\mathbb{C}(z)$. If $r=2$, then the functions (\ref{a21}) are 
algebraically independent apart from the following cases: \\
1) $(b_h)$ is constant and $\alpha_{0,2}=1$; \\
2) $(c_h)$ is constant and $\beta_{0,2}=1$; \\
3) $\alpha_{0,2}=\zeta^2, \beta_{0,2}=\zeta^4$; \\
4) $\alpha_{0,2}=\zeta^4, \beta_{0,2}=\zeta^2$, \\
where $\zeta=e^{\pi i/3}$ as in Lemma 3. In the cases 1) and 2), either 
$\varphi_{0,2}(z)\in\mathbb{C}(z)$, or $\lambda_{0,2}(z)\in\mathbb{C}(z)$, respectively, but all other 
functions in (\ref{a21}) are algebraically independent over $\mathbb{C}(z)$. In the cases 3) and 4), after 
removing the functions $\varphi_{0,2^j}(z)\ ($or $\lambda_{0,2^j}(z)) \ (j\in\mathbb{N})$ from 
(\ref{a21}), the remaining functions are algebraically independent over $\mathbb{C}(z)$}.  \\

\textit{Proof.} We suppose that there exist $L,s\in\mathbb{N}$ such that the functions
\[
\gamma_{\mu,r^j}(z),\ \varphi_{\ell,r^j}(z),\ \lambda_{\ell,r^j}(z)  \quad (\mu=1,\ldots,r-1,\ -L\le\ell\le L,\
j=1,\ldots,s)
\]
are algebraically dependent, and deduce a contradiction.

We shall use Lemma 5, where $(b_{1,h})=(a_h), t_1=r-1,\; (b_{2,h})=(b_h), t_2=2L+1,\;(b_{3,h})=(c_h), 
t_3\!=\!2L+1$. Further, $f_{1,m}(z), f_{2,m}(z), f_{3,m}(z)$ are the functions 
$\gamma_\mu(a,z)\, (1\le\mu<r),\ \varphi_\ell(a,z) \ (-L\le\ell\le L),\ \lambda_\ell(a,z) \ (-L\le\ell\le L)$, 
respectively, and $f_{1,m,r}(z), f_{2,m,r}(z),\\ f_{3,m,r}(z)$ are the functions 
$\gamma_{\mu,r}(z) \ (\mu=1,\ldots,r-1),\ \varphi_{\ell,r}(z) \ (-L\le\ell\le L),\ \lambda_{\ell,r}(z) \ 
(-L\le\ell\le L)$. If $r\ge3$, then a contradiction proving Lemma 6 follows immediately from Lemmas 2, 3 
and 5. The same holds, if $r=2$ and $\alpha_{0,2}\ne1, \beta_{0,2}\ne1$. \\

Assume now $r=2$ and $\alpha_{0,2}=1$ (the case $\beta_{0,2}=1$ being similar). Let $f_{2,1}(a,z)
=\varphi_{0,1}(a,z)$ and $f_{2,1,2}(z)=\varphi_{0,2}(z)$. We may proceed as in the proof of Lemma 5, 
but now we assume that the function $f_{2,1,2}(z)$ is omitted if $(b_h)$ is constant, so $a_{2,1,1}=0$ 
in this case. Now $f_{2,1}(1,z)=\varphi_{0,1}(1,z)\in\mathbb{C}(z)$, and thus (\ref{a17}) yields
\[
G(z)-A_{2,1,0}f_{2,1}(1,z)=\sum_{i=0}^{s!p-1}G_i^*(z)\in\mathbb{C}(z),
\]
where $G_i^*(z)=G_i(z)\, (1 \leq i \leq s!p - 1)$ but $G_0^*(z)=G_0(z)-A_{2,1,0}f_{2,1}(1,z)$. We now define $J^*$ as $J$ in the 
proof of Lemma 5, but this time $A_{2,1,0}$ is replaced by 0 in this definition (thus we may have $0 \in J$ but $0 \notin J^*$). If $J^*$ is not empty, we get a contradiction as in the proof of Lemma 5. Therefore we now deduce that $A_{k,m,i}=0$ for all $(k,m,i)\ne(2,1,0)$. 
Thus, $A_{2,1,0}\ne0$ and, by (\ref{a16}),
\[
\sum_{j=1}^s a_{2,1,j}b_h^{(j)}=A_{2,1,0} \quad (h\in\mathbb{N}_0).
\]
If $(b_h)_{h\ge0}$ is not constant, then the sequences $(b_h^{(j)})_{h\ge0} \ (1\le j \le s)$ and 
$(1)_{h\ge0}$ are linearly dependent contrary to property b) from the beginning of this section. If 
$(b_h)_{h\ge0}$ is constant, then $a_{2,1,1}=0$, and hence $(b_h^{(j)})_{h\ge0}\; (2\le j\le s)$ and 
$(1)_{h\ge0}$ are linearly dependent, contrary to property a). Thus, Lemma 6 holds. \hfill $\Box$  \\ 

\textit{Remark 4.} Similarly, by using Remarks 1 and 3, we see that the results of Lemma 6 hold also 
if $\alpha_\ell$ and $\beta_\ell$ satisfy the conditions of Remark 1. Here the cases 3) and 4) do not occur.

\section{Proof of Theorem 3}

Assume, contrary to Theorem 3, that the numbers (\ref{a6}) are algebraically dependent. Then 
there exist a finite set $\mathcal{R}_0\subset\mathbb{N}\setminus\{1\}$ and a positive integer $L_0$ 
such that the numbers
\[
\Gamma_{\mu,r}(\alpha),\ \Phi_{\ell,r}(\alpha),\ \Lambda_{\ell,r}(\alpha) \quad (-L_0\le\ell\le L_0,\ 
r\in\mathcal{R}_0,\ \mu=1,\ldots, d(r)-1)
\]
are algebraically dependent (in the exceptional cases, the numbers mentioned there are removed). 
We shall prove that this is impossible.

From definition (\ref{a2}) of $\mathcal{D}$, it follows that $\mathbb{N}\setminus\{1\}= 
\{d^j: d\in\mathcal{D}, j\in\mathbb{N}\}$ and $\log d_1/\log d_2\notin\mathbb{Q}$ for all 
distinct $d_1,d_2\in\mathcal{D}$, see \cite[p.105]{NTT}. Thus, there exist a finite subset 
$\mathcal{D}_0\subset\mathcal{D}$ and a positive integer $s$ such that 
$\mathcal{R}_0\subset \{d^j: d\in\mathcal{D}_0, j=1,\ldots,s\}$.

Let $L$ be a positive integer. Assume that, for each $d\in\mathcal{D}_0$, we have 
non-zero algebraic numbers $\alpha_{\ell,d}, \beta_{\ell,d} \ (-L\le\ell\le L)$ as in Lemma 6 
above. If $\alpha$ is an algebraic number with $0<|\alpha|<1$, then we 
may choose an $h_0\in\mathbb{N}$ in such a way that $|\alpha|^{2^{h_0}}<
\min(|\alpha_{\ell,d}|, |\beta_{\ell,d}|: -L\le\ell\le L, d\in\mathcal{D}_0)$. If the functions
\[
\widetilde\gamma_{\mu,d}(z),\ \widetilde\varphi_{\ell,d}(z),\ \widetilde\lambda_{\ell,d}(z)
\]
are defined as in (\ref{a18})-(\ref{a20}) above, but now with summation starting from $h_0$, then all numbers
\begin{equation}\label{a22}
\widetilde\gamma_{\mu,d^j}(\alpha),\ \widetilde\varphi_{\ell,d^j}(\alpha),\ 
\widetilde\lambda_{\ell,d^j}(\alpha) \quad (-L\le\ell\le L,\ d\in\mathcal{D}_0,\ \mu=1,\ldots,d-1,\ j=1,\ldots,s) 
\end{equation}
are defined. For these numbers we have the following \\

\textbf{Lemma 7}. \textit{Let $\alpha_{\ell,d}$ and $\beta_{\ell,d}$ satisfy the conditions of Lemma 6. 
Then the numbers (\ref{a22}) are algebraically independent, unless $2 \in \mathcal{D}_0$ and 
one of the following cases hold: \\
1) $(b_h)$ is constant and $\alpha_{0,2}=1$; \\
2) $(c_h)$ is constant and $\beta_{0,2}=1$; \\ 
3) $\alpha_{0,2}=\zeta^2, \beta_{0,2}=\zeta^4$; \\ 
4) $\alpha_{0,2}=\zeta^4, \beta_{0,2}=\zeta^2$. \\
In the cases 1) and 2), either $\widetilde\varphi_{0,2}(\alpha)$ 
or $\widetilde\lambda_{0,2}(\alpha)$, respectively, is algebraic, but all other numbers in 
(\ref{a22}) are algebraically independent. In the cases 3) and 4), after removing the numbers 
$\widetilde\varphi_{0,2^{j}}(\alpha)$ $($or $\widetilde\lambda_{0,2^{j}}(\alpha)) \ (j=1,\ldots,s)$ 
from (\ref{a22}), the remaining numbers are algebraically independent.} \\  

\textit{Proof.} We apply Lemma 4 to the functions
\[
\widetilde\gamma_{\mu,d^j}(z),\ \widetilde\varphi_{\ell,d^j}(z),\ \widetilde\lambda_{\ell,d^j}(z) \quad 
(-L\le\ell\le L,\ d\in\mathcal{D}_0,\ \mu=1,\ldots,d-1,\ j=1,\ldots,s).
\] 
(In the exceptional cases 1) and 2), either $\widetilde\varphi_{0,2}(z)$ or $\widetilde\lambda_{0,2}(z)$, 
respectively, is removed. In the cases 3) and 4), $\widetilde\varphi_{0,2^{j}}(z) \ 
$(or $\widetilde\lambda_{0,2^{j}}(z)) \ (j=1,\ldots,s)$ are removed.) This gives a system of functional 
equations of the type used in \cite[Lemma 2.1]{NTT}. Further, $\log d_i^{s!p}/\log d_k^{s!p}
\notin\mathbb{Q}$ holds for all distinct $d_i,d_k\in\mathcal{D}_0$. Therefore, by \cite[Lemma 2.1]{NTT}, 
to prove Lemma 7, it suffices to show that, for given $d\in\mathcal{D}_0$, the functions
\[
\widetilde\gamma_{\mu,d^j}(z),\ \widetilde\varphi_{\ell,d^j}(z),\ \widetilde\lambda_{\ell,d^j}(z) \quad 
(-L\le\ell\le L,\ \mu = 1,\ldots,d-1,\ j=1,\ldots,s)
\] 
are algebraically independent over $\mathbb{C}(z)$. Therefore Lemma 6 immediately implies Lemma 7.
\hfill $\Box$  \\

\textit{Remark 5}. By using Remark 4, we similarly get the validity of Lemma 7 if the 
hypotheses of Remark 1 are satisfied. \\

By choosing $L=L_0$ and using Lemma 7, we see that the assumption from the beginning of this 
section leads to a contradiction. This proves Theorem 3.

\section{Proof of Theorems 1 and 2}

It is enough to prove Theorem 2 since Theorem 1 is a consequence of it. We first note that, by 
using the definition of $R_n$ and $S_n$, we get
\begin{equation}\label{a23}
R_{kr^h+\ell}=E_{\ell,r}\gamma_1^{kr^h}(\gamma_1^{-2kr^h}-e_{\ell,r}), \ E_{\ell,r}=
\delta^{kr}g_2\gamma_2^\ell, \ e_{\ell,r}=-\delta^{kr+\ell}g_1g_2^{-1}\gamma_1^{2\ell},
\end{equation}
\begin{equation}\label{a24}
S_{kr^h+\ell}=F_{\ell,r} \gamma_1^{kr^h}(\gamma_1^{-2kr^h}-f_{\ell,r}), \ F_{\ell,r}=
\delta^{kr}h_2\gamma_2^\ell, \ f_{\ell,r}=-\delta^{kr+\ell}h_1h_2^{-1}\gamma_1^{2\ell}.
\end{equation}
Here $e_{\ell,r^j}, f_{\ell,r^j}, E_{\ell,r^j}$ and $F_{\ell,r^j} \ (j\in\mathbb{N})$ do not 
depend on $j$, and $|e_{\ell,r}|$ and $|f_{\ell,r}|$ do not depend on $r$. We first consider the case 
$\Omega=\Delta(\gamma_1/\gamma_2)^{\ell_1}=\Delta\delta^{\ell_1}\gamma_1^{2\ell_1}, 
\left|\Delta\right| = 1, \Delta \neq 1$, which means that $e_{\ell,r}=\Delta f_{\ell_1+\ell,r}, 
|e_{\ell,r}|=|f_{\ell_1+\ell,r}|$ for all $\ell,r$.\\

\textit{Remark 6.} The hypothesis $\Delta\ne1$ in Theorem 2 is needed since otherwise $h_2\gamma_2^{\ell_1}R_{kr^h+\ell}
=g_2S_{kr^h+\ell+\ell_1}$ for all $h,\ell$. If $(b_h) = (c_h) = (1)$, for example, this gives a dependence relation $g_2R_{\ell,r} = h_2\gamma_2^{\ell_1}S_{\ell+\ell_1,r}$. \\  

Assume now that the numbers (\ref{a5}) are algebraically dependent. Then there exist a finite set 
$\mathcal{R}_0\subset\mathbb{N}\setminus\{1\}$ and a positive integer $L_0$ such that the numbers
\[
Q_{\mu,r},\ R_{\ell,r},\ S_{\ell,r} \quad (-L_0 \le\ell\le L_0, \ r\in\mathcal{R}_0, \ \mu=1,\ldots,d(r)-1)
\]
are algebraically dependent (in the exceptional cases we remove the numbers mentioned in 
Theorem 2). As in the last section, there exist a finite subset $\mathcal{D}_0\subset\mathcal{D}$ 
and a positive integer $s$ such that $\mathcal{R}_0\subset\{d^j\!: d\in\mathcal{D}_0, j=1,\ldots,s\}$.

Let $\alpha_{\ell,r}:=e_{\ell_0+\ell,r},\ \beta_{\ell,r}:=f_{\ell_1+\ell_0+\ell,r}$ for $\ell\in\mathbb{Z}$ if 
there exists some $\ell_0$ with $|e_{\ell_0,r}|=1$; otherwise we choose $\ell_0=0$. Let 
$L=L_0+|\ell_1|+|\ell_0|$, and choose $h_0$ such that
\[
|\gamma_1|^{-k2^{h_0}}<\min(|\alpha_{-L,r}|, |\beta_{-L,r}|).
\]
According to Lemma 7, the numbers
\[
\widetilde\gamma_{\mu,d^j}(\gamma_1^{-k}),\ \widetilde\varphi_{\ell,d^j}(\gamma_1^{-k}),\ 
\widetilde\lambda_{\ell,d^j}(\gamma_1^{-k}) \quad (-L\le\ell\le L,\: d\in\mathcal{D}_0,\: 
\mu=1,\ldots,d-1,\: j=1,\ldots,s)
\] 
are algebraically independent. (In the exceptional cases 1), 2), and 3,4), the numbers 
$\widetilde\varphi_{0,2}(\gamma_1^{-k}), \\ 
\widetilde\lambda_{0,2}(\gamma_1^{-k})$ and $\widetilde\varphi_{0,2^j}(\gamma_1^{-k}) \ 
(j=1,\ldots,s)$, respectively, are removed.) Since, by (\ref{a23}) and (\ref{a24}), for any $r=d^j$,
\[
Q_{\mu,r}-\sum_{h=0}^{h_0-1}\frac{a_h}{\gamma_1^{k\mu r^h}}=\widetilde\gamma_{1,r}(\gamma_1^{-k}),
\]
\[
R_{\ell,r}-{\sum_{h=0}^{h_0-1}}\strut'\frac{b_h}{R_{kr^h+\ell}}=
E_{\ell,r}\widetilde\varphi_{\ell-\ell_0,r}(\gamma_1^{-k}),
\]
\[
S_{\ell,r}-{\sum_{h=0}^{h_0-1}}\strut'\frac{c_h}{S_{kr^h+\ell}}=
F_{\ell,r}\widetilde\lambda_{\ell-\ell_1-\ell_0,r}(\gamma_1^{-k}),
\]
we get a contradiction to the assumption on algebraic dependence above.
Furthermore, note that the condition $R_{\ell_0}=0$ \ (or $S_{\ell_1+\ell_0}=0$) is equivalent 
to $\alpha_{0,2}=e_{\ell_0,2}=1$ \ (or $\beta_{0,2}=f_{\ell_1+\ell_0,2}=1$), and the conditions in 3,4) 
mean that $\alpha_{0,2}=e^{2\pi i/3}$ or $e^{4\pi i/3}$ and $\beta_{0,2}=\alpha_{0,2}^2$.

If $|\Omega|\notin|\gamma_1/\gamma_2|^{\mathbb{Z}}$, then $|e_{\ell_1,r}|\ne|f_{\ell_2,r}|$ for 
all $\ell_1,\ell_2$. Thus, the above deduction works also in this case, by Remark 5. This proves 
Theorem 2. \hfill $\Box$ 

\vskip0.8cm
\begin{small}
\noindent
Authors' addresses:   \\
\vskip0.1cm
\noindent
Peter Bundschuh \hfill Keijo V\"a\"an\"anen  \\
Mathematisches Institut  \hfill  Department of Mathematical Sciences \\ 
Universit\"at zu K\"oln  \hfill  University of Oulu \\
Weyertal 86-90  \hfill  P. O. Box 3000  \\
50931 K\"oln, Germany  \hfill   90014 Oulu, Finland \\
E-mail: pb@math.uni-koeln.de  \hfill   E-mail: kvaanane@sun3.oulu.fi \\
\end{small}

\end{document}